\documentclass[french,10pt]{article}
\usepackage[T1]{fontenc}
\usepackage[french]{babel}
\usepackage{amsfonts}\usepackage{amssymb}
\usepackage{euscript}
\usepackage{diagrams}
\usepackage{mathrsfs} 
\usepackage{amsmath}
\newtheorem{proposition}{Proposition}[section]
\newtheorem{df}{D\'efinition}[section]
\newtheorem{thm}{Th\'eor\`eme}[section]
\newtheorem{lemme}{Lemme}[section]
\newtheorem{cor}{Corollaire}[section]

\def\h{\hbox{$/ \kern-.3em /$}}
\newcommand{\sq}[0]{ {\hfill$\square$} }
\newcommand{\ssq}[0]{ {\hfill$\square\ \square$} }
\newcommand{\preuve}[0]{\noindent{{\em Preuve}}}
\newcommand{\preuvethm}[0]{\vspace{0.2 cm}\noindent{\em{Preuve du th\'eor\`eme}}}

\title{ Sur la cat\'egorie $\mathcal{D}^G(X)$ pour l'action  
d'un groupe fini \\avec quotient lisse}
\date{Avril 2002}\author{Sophie T\'erouanne}
\begin{document}

\maketitle
\begin{center}
Pr\'epublication de l'Institut Fourier \FrenchEnumerate{n} 569 (2002)
\\{\em http://www-fourier.ujf-grenoble.fr/prepublications.html}
\end{center}
\begin{center}
\parbox[c]{10.5 cm}{
{\small
\begin{center}{\bf Abstract}
\end{center}
Let $X$ be a smooth projective scheme over an algebraically closed field $k$  
of characteristic $0$. Let $G$ be a finite group acting on $X$.
If $X/G$ is smooth, then $X/G\stackrel{\rm{can}}{\cong} X\h G$ but the theorem of \cite{BKR:01}
does not a priori apply because the dualizing sheaf of 
$X$ is not trivial as a $G$-sheaf.
We compare the categories $\mathcal{D}^G(X)$ and $\mathcal{D}(X/G)$ 
in this case by giving the image of the Fourier-Mukai functor
$\begin{diagram} \mathcal{D}(X/G) \simeq \mathcal{D}(X\h G) &\rTo &\mathcal{D}^G(X) \end{diagram}$
defined in \cite{BKR:01}. The result generalizes a theorem of L\o nsted on
$G$-equivariant $K$-theory in the case of curves (\cite{Lonsted:83})
}}
\end{center}

%%%%%%%%%%%%%%%%%%%%%%%%%%%%%%%%%%%%%%%%%%%%%%%%%%%%%%%%%%%%%%%%%%%%%%%%%%%%%%%

%%%%%%%%%%%%%%%%%%%%%%%%%%%%%%%%%%%%%%%%%%%%%%%%%%%%%%%%%%%%%%%%%%%%%%%%%%%%%%

\section*{Introduction}

Pour un groupe fini $G$ agissant sur un sch\'ema lisse $X$ donn\'e, 
la cat\'egorie d\'eriv\'ee $G$-\'equivariante de $X$
est par d\'efinition la 
cat\'egorie d\'eriv\'ee born\'ee des $G$-faisceaux
quasicoh\'erents sur $X$, 
que l'on a restreinte  aux complexes \`a homologie coh\'erente.
On la notera $\mathcal{D}^G(X)$.
%Cette cat\'egorie contient des informations sur le sch\'ema $X$, l'action 
%de $G$ sur $X$.
Notons $X\h G$ la composante irr\'eductible du $G$-sch\'ema de Hilbert
 $G{\rm -Hilb}(X)$
contenant les $G$-orbites libres.
Sous une  hypoth\`ese sur la  dimension de 
$X\h G \times_{X/G} X\h G$, et lorsque le faisceau
dualisant de $X$ est trivial en tant que $G$-faisceau,
Bridgeland, King et Reid ont d\'emontr\'e (\cite{BKR:01}) 
que  $X\h G$
est une r\'esolution cr\'epante du quotient $X/G$ et que
la cat\'egorie $\mathcal{D}^G(X)$ est \'equivalente \`a la cat\'egorie d\'eriv\'ee 
$\mathcal{D}(X\h G)$. 

Il semble alors l\'egitime de se poser la question suivante : 
lorsque le quotient $X/G$ est lisse, et qu'il n'est donc pas n\'ecessaire
de le d\'esingulariser, y-a-t-il toujours \'equivalence entre la cat\'egorie d\'eriv\'ee
$G$-\'equivariante de $X$,  et la cat\'egorie d\'eriv\'ee du 
quotient $\mathcal{D}(X/G)$  ?

Dans cet article,  $X$ est un sch\'ema projectif lisse sur un corps $k$ 
alg\'ebriquement clos de 
caract\'eristique nulle. On suppose qu'un groupe fini $G$ 
agit sur $X$ de façon \`a ce que le quotient $X/G$ soit lisse. 
Par exemple, $X$ est le produit  $E^n$ où $E$ est une courbe lisse et 
$G = \mathfrak{S}^n$ agit par permutation des composantes.
On note
$\begin{diagram}[size=2em]\pi: X & \rTo &X/G\end{diagram}$  le morphisme quotient. 
Les foncteurs adjoints suivant apparaissent alors naturellement. 
$$\begin{diagram}\pi_*^G : \mathcal{D}^G(X) &\rTo& 
\mathcal{D}(X/G)\end{diagram}$$
$$\begin{diagram}\pi^* : \mathcal{D}(X/G)&\rTo& 
\mathcal{D}^G(X)\end{diagram}.$$ 
Ces foncteurs sont exacts entre les cat\'egories de faisceaux car, par lissit\'e
de $X/G$, 
$\pi$ est plat, et on est
sur un corps $k$ dont la caract\'eristique ne divise pas l'ordre de $G$.
Ainsi, ils induisent naturellement des foncteurs entre les cat\'egories
d\'eriv\'ees.

En consid\'erant $X$ comme une famille plate de $G$-clusters de $X$
index\'ee par  $X/G$,  on d\'emontre ais\'ement que 
$X/G\stackrel{\rm can}{\cong} X\h G$ avec $X$ pour famille universelle. 
Ainsi, le foncteur de Fourier-Mukai 
$\begin{diagram}[size=2em] \mathcal{D}(X/G) & \rTo & \mathcal{D}^G(X)\end{diagram}$ 
d\'efini dans \cite{BKR:01} est le foncteur $\pi^*$.
Le fait que le quotient $X/G$ soit lisse implique que localement (en tant que
sous groupe du groupe d'endomorphisme de l'espace tangent), $G$ est engendr\'e par
des pseudor\'eflexions (\cite{Lie}). En particulier, le faisceau canonique de $X$ n'est 
pas $G$-trivial a priori, de sorte que l'on ne peut pas appliquer le th\'eor\`eme
donnant une \'equivalence de cat\'egories.             
Le r\'esultat de cet article est le th\'eor\`eme suivant  qui 
caract\'erise l'image de $\pi^*$ dans $\mathcal{D}^G(X)$.

\vspace{0.2 cm}\noindent{\bf Th\'eor\`eme (\ref{thm_critere_general})\ }{\it
Un objet  $E\in\mathcal{D}^G(X)$ est l'image par $\pi^*$ d'un objet
de $\mathcal{D}(X/G)$
 si et seulement 
si pour tout $H$ sous-groupe de $G$, 
la restriction 
de $E$ au sous-sch\'ema $\overline{\Delta}_H$ (adh\'erence
des points de stabilisateur $H$) 
appartient \`a la composante correspondant  \`a la repr\'esentation 
triviale dans la d\'ecomposition de la cat\'egorie d\'eriv\'ee 
$ \mathcal{D}^{H}(\overline{\Delta}_H) = 
\oplus_{\rho\in{\rm Irr}(H)}\mathcal{D}^\rho(\overline{\Delta}_H)$.
}

\vspace{0.2 cm}

Dans la premi\`ere partie, nous d\'efinissons la cat\'egorie 
$\mathcal{D}^G(X)$ et nous d\'emontrons qu'elle est \'equivalente
\`a la sous-cat\'egorie form\'ee des $G$-faisceaux coh\'erents localement libres. 
Nous \'etudions les cas extrêmes où $G$ agit soit trivialement, soit 
librement sur $X$. 
Puis, en suivant le mod\`ele donn\'e
par Be{\u\i}linson dans \cite{Bej:79} dans le cas non \'equivariant,
nous donnons  
une description de $\mathcal{D}^G(\mathbb{P}^n)$, 
où $\mathbb{P}^n$ est l'espace projectif de dimension $n$ sur $k$.

Dans la deuxi\`eme partie, nous  d\'emontrons 
le th\'eor\`eme~\ref{thm_critere_general} 
apr\`es avoir montr\'e comment r\'eduire le probl\`eme \`a la recherche
d'un crit\`ere de descente pour les  
$G$-faisceaux localement libres sur $X$.

%%%%%%%%%%%%%%%%%%%%%%%%%%%%%%%%%%%%%%%%%%%%%%%%%%%%%%%%%%%%%%%%%%%%%%%%%%%%%%%
%%%%%%%%%%%%%%%%%%%%%%%%%%%%%%%%%%%%%%%%%%%%%%%%%%%%%%%%%%%%%%%%%%%%%%%%%%%%%%%

\section{La cat\'egorie d\'eriv\'ee born\'ee $G$-\'equivariante}
%%%%%%%%%%%%%%%%%%%%%%%%%%%%%%%%%%%%%%%%%%%%%%%%%%%%%%%%%%%%%%%%%%%%%%%%%%%%%%%
\subsection{D\'efinition}

Rappelons que si $X$ est un sch\'ema et $G$ un 
groupe fini agissant sur $X$, 
un faisceau $F$ coh\'erent est dit $G$-lin\'earis\'e
si il existe une collection d'isomorphismes :  
$\begin{diagram}[size=2em]\phi _g :F&\rTo& g_* F,\  g\in G\end{diagram}$,  telle que 
$\phi_{gh}^{-1} = h_*\phi_g^{-1}\phi_h^{-1}$.
On appelle morphisme de $G$-faisceaux un morphisme 
de faisceaux entre deux faisceaux $G$-lin\'earis\'es
qui commute avec les morphismes $\phi_g,$ $g\in G$.
Alors on d\'efinit la cat\'egorie $G{\rm{-Coh}}(X)$ dont 
les objets sont les faisceaux coh\'erents $G$-lin\'earis\'es et les
morphismes sont les morphismes de  $G$-faisceaux. 
Cette cat\'egorie est ab\'elienne et elle a assez d'injectifs 
(cf {\cite{Grothendieck:57}}).

\begin{lemme}\label{dhn}
Si $X$ est un $G$-sch\'ema lisse de dimension $n$, alors  la cat\'egorie 
$G{\rm{-Coh}}(X)$ est de dimension homologique inf\'erieure \`a $n$.
C'est \`a dire, pour tout couple de 
$G$-faisceaux coh\'erents $(A,B)$ et pour tout entier $s\ge n$,  ${\rm{Ext}}^s_G(A,B) = 0$.
\end{lemme}

\preuve .
Notons  ${\rm I}^G$ le foncteur 
qui \`a un $G$-module associe sa partie $G$-invariante, 
et $\Gamma$ le foncteur qui \`a un $G$-faisceau associe
ses sections globales.
Notons $H_G^*$ la cohomologie $G$ \'equivariante. 
Alors, le foncteur ${\rm{Ext}}_G$ est le foncteur d\'eriv\'e 
de ${\rm I}^G\circ\Gamma\circ{\mathscr{H}}om$, et 
pour tout couple $(A,B)$ de $G$-faisceaux coh\'erents,
la suite spectrale de ter\-me initial 
${\rm E}^{p,q}_2 = H_G^p(X,     \mathscr{E}xt^q (A,B)$
%{\PolicePourFophie Ext}^q(A,B))$ 
converge vers 
${\rm{Ext}}^{p+q}_G(A,B)$.
On va d\'emontrer que pour tout couple d'entiers 
$(p,q)$ tel que $p+q> n$, tous les 
termes de cette suite spectrale sont nuls.

Fixons provisoirement $q$, et 
notons $d$ la dimension du support du faisceau $\mathscr{E}xt^q(A,B)$
%${\PolicePourFophie Ext}^q(A,B)$.
La caract\'eristique de $k$ ne divisant pas l'ordre de $G$, 
${\rm I}^G$ est un foncteur exact sur la cat\'egorie des 
$k$-espaces vectoriels, de sorte que l'on a 
$ H_G^p(X,\mathscr{E}xt^q(A,B)          % ${\PolicePourFophie Ext} $\!^q(A,B)
) = 
{\rm I}^G( H^p(X,\mathscr{E}xt^q(A,B) % ${\PolicePourFophie Ext} $\!^q(A,B)
))$.
Par ailleurs, si 
$\begin{diagram}[size=2em]
j : Z:= {\rm {Supp}}(\mathscr{E}xt^q(A,B) %${\PolicePourFophie Ext} $\!^q(A,B)
)& \rTo & X\end{diagram}$
est l'inclusion du support de $\mathscr{E}xt^q(A,B)$ %{\PolicePourFophie Ext} $\!^q(A,B)$ 
dans $X$, alors
on a (\cite[p.209]{Hartshorne:66}) 
$ H^p(X,j_\star(\mathscr{E}xt^q(A,B)_{|Z}))
%${\PolicePourFophie Ext} $\!^q(A,B)_{|Z})) 
= H^p(Z,\mathscr{E}xt^q(A,B)_{|Z} % ${\PolicePourFophie Ext} $\!^q(A,B)_{|Z}
)
=  H^p(Z,\mathscr{E}xt^q(A,B) %${\PolicePourFophie Ext} $\!^q(A,B)
)$
donc par th\'eor\`eme d'annulation de Grothendieck, $H^p(X,\mathscr{E}xt^q(A,B)
% ${\PolicePourFophie Ext} $\!^q_G(A,B)
)=0$ si $p>d$.

Supposons maintenant que %{\PolicePourFophie Ext} $\!^q(A,B)$ 
$\mathscr{E}xt^q(A,B) $ est  non nul et 
majorons la dimension $d$ de $Z$.
Soit $x$ un point g\'en\'erique de $Z$. Alors ${\mathscr{O}}_{X,x}$ 
est un anneau r\'egulier de dimension inf\'erieure \`a $n-d$.
Comme $x$ est dans le support de %{\PolicePourFophie Ext} $\!^q(A,B)$ 
$\mathscr{E}xt^q(A,B) $ par hypoth\`ese, 
la fibre $\mathscr{E}xt^q(A,B)_x  % ${\PolicePourFophie Ext} $\!^q(A,B)_x 
= {\rm Ext}^q_{\mathscr{O}_{X,x}}(A_x, B_x)$
est un $\mathscr{O}_{X,x}$-module non nul. Or la dimension 
homologique de la cat\'egorie des  $\mathscr{O}_{X,x}$-module
est inf\'erieure    \`a la dimension de  $\mathscr{O}_{X,x}$, 
donc la non-nullit\'e de ${\rm Ext}^q_{\mathscr{O}_{X,x}}(A_x, B_x)$
implique $q\leq n-d$.

Finalement, pour que ${\rm E}^{p,q}_2 = H_G^p(X, \mathscr{E}xt^q(A,B)  
% ${\PolicePourFophie Ext} $\!^q(A,B)
)$ soit 
non nul, il est n\'ecessaire d'avoir $p+q\leq n$.
En particulier, pour tout entier $s>n$, ${\rm{Ext}}^s_G(A,B) = 0$.
\sq
%%%%%%%%%%%%%%%%%%%%%%%%%%%%%%%%%%%%%%%%%%%%%%%%%%%%%%%%%%%%%%%%%%%%%%%%%%%%%%%
%%%%%%%%%%%%%%%%%%%%%%%%%%%%%%%%%%%%%%%%%%%%%%%%%%%%%%%%%%%%%%%%%%%%%%%%%%%%%%%
\subsection{Une \'equivalence de cat\'egories}\label{red1}

Pour d\'ecrire la cat\'egorie d\'eriv\'ee $G$-\'equivariante, on va voir 
qu'il suffit de connaître la sous-cat\'egorie pleine 
dont les objets sont des complexes de faisceaux 
coh\'erents localement libres.

\begin{thm}\label{equivalence_cat_der_sous_cat_loclibre}
Soit $X$ un sch\'ema projectif lis\-se et soit $G$ un groupe fini agissant sur $X$.
Soit $\mathscr{L}\subset G{\rm -Coh}(X))$ la sous-cat\'egorie pleine des
$G$-faisceaux  coh\'erents localement libres sur $X$ et 
$\mathscr{M}\subset \mathcal{D}^G(X)$ la  sous-cat\'egorie pleine des complexes
de faisceaux  dans $\mathscr{L}$.
On a une \'equivalence de cat\'egories $$\mathscr{M}\simeq \mathcal{D}^G(X).$$
\end{thm}

Comme on a impos\'e dans la d\'efinition que $\mathscr{M} $ soit une 
sous-cat\'egorie pleine de  
$\mathcal{D}^G(X)$, il suffit de d\'emontrer que tout objet de  
$\mathcal{D}^G(X)$
admet un quasi-isomorphisme 
%de la cat\'egorie d\'eriv\'ee $\mathcal{D}^G(X)$ 
vers un objet de $\mathscr{M}$. 
Remarquons que $\mathscr{M}$ n'est pas la cat\'egorie d\'eriv\'ee de $\mathscr{L}$.
En effet, en imposant que la sous-cat\'egorie soit pleine, on a ajout\'e 
des  morphismes par rapport \`a ceux de $\mathcal{D}(\mathscr{L})$.
Soit donc $F\in {\rm Ob}(\mathcal{D}^G(X))$. On va d\'emontrer 
par r\'ecurrence sur 
l'amplitude homologique  de $F$ qu'il est quasi-isomorphe \`a un 
objet de $\mathscr{M}$.

Si $F$ est d'amplitude nulle, alors il admet au plus un 
faisceau  d'homologie non nul. Et si celui-ci est en degr\'e $i$, alors 
$F\cong H^i(F)[-i]$.   
Il suffit donc de d\'emontrer que le faisceau $H^i(F)$ admet une 
r\'esolution libre finie. C'est le lemme classique suivant.
%%%%%%%%%%%%%%%%%%%%%%%%%%%%%%%%%%%%%%%%%%%%%%%%%%%%%%%%%%%%%%%%%%%
% SUIVANT QUI ??????-> trouver la reference.
%%%%%%%%%%%%%%%%%%%%%%%%%%%%%%%%%%%%%%%%%%%%%%%%%%%%%%%%%%%%%%%%%%%%%%%%%%%%%%%

\begin{lemme}
Soit $X$ un sch\'ema projectif sur laquelle agit $G$ un groupe fini.
Tout $G$-faisceau coh\'erent admet une r\'esolution finie par des 
$G$-faisceaux localement
libre de types  finis.
\end{lemme}
\preuve .
Pour commencer, on remarque que $X$ admet un faisceau tr\`es ample 
$G$-lin\'earis\'e. 
En effet, on peut plonger $X$ dans un espace projectif $\mathbb{P}(V)$ 
de telle 
façon que 
$G$ agisse lin\'eairement sur $V$. Soit $L$ un faisceau tr\`es ample sur
$\mathbb{P}(V)$. A priori, $L$ n'est pas $G$-lin\'earis\'e. En revanche, 
le produit 
$\otimes_{g\in G}g^*(L)$ est ample et admet une $G$-lin\'earisation naturelle.
 Il existe une puissance de ce  faisceau $M = L^{\otimes n}$ qui est 
tr\`es ample. Le pull-back de $M$ sur $X$ (encore not\'e $M$) est bien un 
$G$-faisceau 
tr\`es ample sur $X$.

Soit $F$ un $G$-faisceau  coh\'erent sur $X$. Par th\'eor\`eme, il existe
un entier $n$ tel que $F\otimes M^{\otimes n}$ est engendr\'e pas ses sections 
globales.
Alors $G$  agit naturellement sur $\Gamma(X, F\otimes M^{\otimes n})$ et  
le morphisme de $G$-faisceaux : 
$$\begin{diagram}
\mathscr{O}_X\otimes \Gamma(X, F\otimes M^{\otimes n})&
                 \rTo &F\otimes M^{\otimes n}\\ a\otimes s&\rTo &as
\end{diagram}
$$ est surjectif.
D'où la surjection $$\begin{diagram}
\Gamma(X, F\otimes M^{\otimes n})\otimes M^{\otimes (-n)}&\rTo& F\end{diagram}.$$
En r\'ep\'etant autant de fois que n\'ecessaire
ce proc\'ed\'e en remplaçant $F$ par le noyau de ce morphisme, 
on obtient une r\'esolution de $F$ par des $G$-faisceaux 
coh\'erents localement libres.
Or d'apr\`es le lemme~\ref{dhn},  la cat\'egorie $G{\rm -Coh}(X)$ est de dimension 
homologique finie, donc cette it\'eration s'arête et  
 $F$ admet une r\'esolution libre finie par des objets de   $\mathscr{L}$.
\sq

%%%%%%%%%%%%%%%%%%%%%%%%%%%%%%%%%%%%%%%%%%%%%%%%%%%%%%%%%%%%%%%%%%%%%%%%%%%%%%%

\preuvethm~\ref{equivalence_cat_der_sous_cat_loclibre}.
Le lemme initialise notre r\'ecurrence sur l'amplitude homologique de $F$.

Supposons que pour tout complexe d'amplitude $k$, il existe un complexe 
d'objets de  $\mathscr{L}$ qui lui est quasi-isomorphe,
 et d\'emontrons le pour un complexe d'amplitude $k+1$.
Soit $F$ un complexe d'amplitude $k+1$. On peut supposer que les faisceaux 
d'homologie non nuls se trouvent en les degr\'es $i$ tels que  $0\leq i\leq k$. 
On consid\`ere le complexe tronqu\'e $\tau^{\ge 1}F$. Il est d'amplitude  
$k-1$, donc
par hypoth\`ese de r\'ecurrence, il existe un complexe $L\in  \mathscr{M}$ et 
un morphisme de complexe  
$\begin{diagram}[size=2em]L&\rTo&\tau^{\geq 1} F\end{diagram}$ 
qui soit un isomorphisme  en cohomologie.
Ainsi on a le diagramme ci-apr\`es où pour tout $i\geq 2$, $\alpha ^i$ est 
un isomorphisme en cohomologie, 
et $\alpha^1(\ker(d^1)) = \ker(\partial^{1})$.
$$\begin{diagram}
  \rTo & F^0 &\rTo^{\partial^0} &F^1 &\rTo^{\partial^1} &\cdots &
            \rTo{\partial^{k-1}}&F^k &\rTo^{\partial^{k}} &F^{k+1}&\rTo \\
             &     &         & \uTo^{\alpha^1}& &       &         &
            \uTo^{\alpha^k}& & \uTo^0 \\
                &     &            &L^1&\rTo^{d^1}& \cdots&\rTo{d^{k-1}}&  
              L^k &\rTo & 0 &\rTo \\
\end{diagram}
$$
Pour commencer, compl\'etons ce morphisme de complexes en construisant: 

$\begin{diagram} F^0 \\ \uTo^{\alpha^0}\\ 
L^0 &\rTo^{d^0}&L^1 \end{diagram}$
tel que :$\left\{\begin{array}{l} \alpha^1 { \rm \ est \ un \ isomorphisme\  
en\ cohomologie,} \\   \alpha^0(\ker(d^0)) = \ker(\partial^{0})= H^0(F).
\end{array}\right.$

Soient $p_1$ et $p_2$ les projections de $F^0\times \ker(d^1)$ 
sur $F^0$ et $\ker(d^1)$. Notons $ F^0\times^{F^1}\ker(d^1)$
le noyau de l'application $\partial^0 p_1 - \alpha^1p_2$.
C'est un sous-$G$-faisceau coh\'erent de $F^0\times \ker(d^1)$.

\vspace{0.1 cm}
\noindent
$
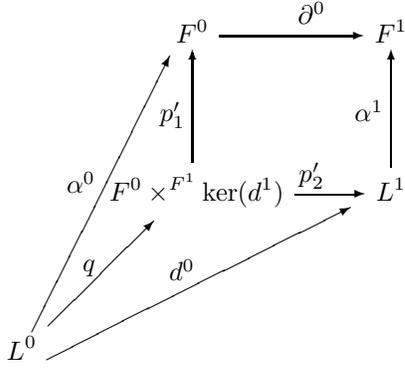
\begin{diagram} 
 &&      F^0  & \rTo^{\partial^0} & F^1\\ 
&\ruTo(2,4)^{\alpha^0}&\uTo^{p_1'}               &           &\uTo^{\alpha^1}\\ 
 &                    & \ F^0\times^{F^1}\ker(d^1)&\rTo^{p'_2} & L^1 \\
&\ruTo^{q}&& \ruTo(4.85,2.25)^{d^0}\\
L^0 & \\
\end{diagram}
$\hspace{0.2 cm}
\parbox[r]{9.5 cm}
{ Notons $p_i'$ la restriction de $p_i$  
\`a \mbox{$F^0\times^{F^1}\ker(d^1)$}. D'apr\`es le lemme, il existe un
$G$-faisceau coh\'erent localement libre $L^0$ et une application  
$q$
surjective.

Posons $\alpha^0 = p_1'\circ q$ et  $d^0 = p_2'\circ q$.
Alors 
${\rm im}(d^0) = {\rm im}(p'_2)$ et un calcul donne 
${\rm im}(p_2')_x = {\rm ker}(d^1)_x\cap ({\alpha}^1)^{-1}({\rm{im}}(\partial^0)_x)
%\{m\in (\ker d^1)_x {\mbox{ tel que }} \alpha^1(m)\in {\rm im}(\partial ^0)_x \}
                 = {\rm ker}(d^1)_x\cap \ker(\bar{\alpha}^1)_x$ 
avec $\bar{\alpha}^1 = \pi\circ\alpha^1$,
où $\pi$ est le quotient de $F^1$ par ${\rm im}(\partial^0)$.
Ainsi, ${\rm im}(d^0) = \ker(\bar{\alpha}^1_{| {\rm ker}d^1})$, et par 
propri\'et\'e universelle 
du noyau, cela implique que $\alpha^1$ est un isomorphisme en cohomologie, 
car $\bar{\alpha}^1_{| {\rm ker}d^1}$ 
induit un isomorphisme 
$$\ker (d^1 ) / {\rm im } (d^0)\cong 
                           \ker(\partial^1) /{ \rm im}(\partial ^0).    $$
}

Puis $\alpha^{0}(\ker (d^0)) = p'_1\circ q(\ker (p'_2\circ q) ) = 
                                                      p'_1(\ker( p'_2))$, 
car $q$ est surjective. En se plaçant sur une fibre, on a 
$p'_1(\ker p'_2)_x = p'_1(F^0_x\times^{F^1_x} \{0\})
                                     = \ker(\partial ^0)_x,  $
d'où la deuxi\`eme condition \`a v\'erifier : 
$$\alpha^0( \ker (d^0) )  = \ker(\partial^0).$$

Alors pour avoir un quasi-isomor\-phisme de complexe, il suffit de 
cons\-truire une r\'esolution 
libre finie 
$\begin{diagram}L^{-k}&\rTo &  \ldots &L^{-1}&\rTo\end{diagram}$ 
du noyau  de $\alpha^0$. En effet, cela implique que
$\alpha^0$ est un isomorphisme en cohomologie puis en posant  
$\alpha^i = 0$ pour tout $i<0$, 
on en d\'eduit que $\begin{diagram}[size=2em]\alpha : L &\rTo & F\end{diagram}$ est un quasi-isomorphisme.
D'apr\`es le lemme cette r\'esolution existe, et elle est finie.

On a donc bien construit un complexe $L\in \mathscr{M}$ qui est 
quasi-isomorphe \`a $F$, 
d'où l'\'equivalence de cat\'egorie.
\sq

%%%%%%%%%%%%%%%%%%%%%%%%%%%%%%%%%%%%%%%%%%%%%%%%%%%%%%%%%%%%%%%%%%%%%%%%%%%%%%%
%%%%%%%%%%%%%%%%%%%%%%%%%%%%%%%%%%%%%%%%%%%%%%%%%%%%%%%%%%%%%%%%%%%%%%%%%%%%%%%

\subsection{Cas d'une action libre }\label{actionlibre}
Dans le cas où l'action de $G$ est libre sur $X$, si $X$ est lisse alors
le quotient $X/G$ est lisse et  on
trouve dans  \cite[p.70]{Mumford:74} la d\'emonstration de
la proposition suivante.
\begin{proposition}\label{libre}
Lorsque $G$ agit librement sur un sch\'ema $X$, on a une \'equivalence 
de cat\'egorie :
 $$ G{\rm{-Coh}}(X)\simeq{\rm{Coh}} (X/G).$$
\end{proposition}

On en d\'eduit alors 
\begin{cor}
Si $G$ agit librement sur un sch\'ema   lisse $X$, le foncteur 
$$\begin{diagram}
\pi^* : \mathcal{D}(X/G) &\rTo& \mathcal{D}^G(X)
\end{diagram}$$
est une \'equivalence de cat\'egorie.
\end{cor}

\preuve . Comme dans le th\'eor\`eme~\ref{equivalence_cat_der_sous_cat_loclibre} 
consid\'erons un triangle du type : 
$$\begin{diagram}
\tau^{\ge {1}} F[-1]&\rTo& H^0(F)&\rTo &F&\rTo & 
                                        \tau^{\ge 1} F\end{diagram},$$ 
et on 
raisonne par r\'ecurrence sur l'amplitude.
On se ram\`ene ainsi \`a d\'emontrer le  lemme suivant :
\begin{lemme}
Soient $E$ un faisceau coh\'erent sur $X/G$  et $F$ un  objet
de $\mathcal{D}(X/G)$. On peut consid\'erer $E$ comme \'etant un complexe
concentr\'e en degr\'e $0$. Si 
$\begin{diagram}[size=2em] \alpha : \pi^*F&\rTo&\pi^* E \end{diagram}$ 
est un morphisme dans $\mathcal{D}^G(X)$, 
alors $\alpha = \pi^*\beta$
pour un certain morphisme  $\beta\in \mathcal{D}(X/G)$.
\end{lemme}
\preuve .
D'apr\`es le th\'eor\`eme~\ref{equivalence_cat_der_sous_cat_loclibre} on 
peut remplacer $F$ par un complexe $L$ 
de faisceaux localement libres et par exactitude de $\pi^*$,
$\pi^*L$ est alors une r\'esolution de $\pi^* F$.
Ainsi, on a 
$$\begin{array}{rcl}
H^i(
%{\mbox{\PolicePourFophie Hom } } 
\mathscr{H}om(\pi^*F,\pi^*E))& = & H^i(%{\mbox{\PolicePourFophie Hom }}
\mathscr{H}om  (\pi^*L,\pi^*E))\\
                                  & = &  %{\PolicePourFophie Ext} 
\mathscr{E}xt^i (\pi^*L,\pi^*E)\\
                                  & = & 0, {\ \rm pour \ i\neq 0, }
\end{array}
$$
car $\pi^* L$ est localement libre. On est donc ramen\'e \`a des morphismes 
de faisceaux 
et on peut utiliser la proposition \ref{libre} : 
$$\begin{array}{rcl}
{\rm Hom} (\pi^*F,\pi^*E)& = &{\rm Hom}_{G{\rm -Coh}(X)} (\pi^*L^0,\pi^*E)\\
                        & \cong &  {\rm Hom}_{{\rm Coh}(X/G)} (L^0,E)\\
                         & \cong & { \rm Hom}_{\mathcal{D}(X/G)} (L,E)\\
                         & \cong & { \rm Hom}_{\mathcal{D}(X/G)} (F,E).\\ 
\end{array}
$$
\ssq

%%%%%%%%%%%%%%%%%%%%%%%%%%%%%%%%%%%%%%%%%%%%%%%%%%%%%%%%%%%%%%%%%%%%%%%%%%%%%%%
%%%%%%%%%%%%%%%%%%%%%%%%%%%%%%%%%%%%%%%%%%%%%%%%%%%%%%%%%%%%%%%%%%%%%%%%%%%%%%%

\subsection{Cas de l'action triviale}\label{decomposition}

On a trait\'e le cas de l'action libre. 
L'autre cas simple et extr\`eme est celui où $G$ agit trivialement sur le 
sch\'ema $X$. 
Cela ne nous empêche pas de d\'efinir la cat\'egorie  des $G$-faisceaux, 
ni a fortiori la cat\'egorie d\'eriv\'ee $\mathcal{D}^G(X)$.
Nous allons voir que dans ce cas, cette cat\'egorie se d\'ecompose en le 
sens suivant :

\begin{df}
On dit qu'une cat\'egorie additive $\mathscr{A}$ se d\'ecompose en 
$\oplus\mathscr{A}_i$  si pour tout $A\in{\rm Ob}(\mathscr{A})$, il existe
$A_i\in \mathscr{A}_i$ tels que 
$A = \oplus A_i$ et ${\rm Hom}_\mathscr{A}(A_i, A_j) = 0$ si $i\neq j$.
\end{df}

\begin{lemme}\label{dec A->dec D(A)}
Soit $\mathscr{A}$ une cat\'egorie ab\'elienne qui se d\'ecompose en 
$\mathscr{A} = \oplus_{i\in I}\mathscr{A}_i$, où $I$ est fini.
Supposons que $\mathscr{A}$  et les sous-cat\'egories  $\mathscr{A}_i$ 
admettent  assez d'injectif et que les injectifs 
des sous-cat\'egories soient des injectifs de $\mathscr{A}$.
Alors la cat\'egorie d\'eriv\'ee $\mathcal{D}(\mathscr{A})$ se d\'ecompose en  
$\oplus_{i\in I}\mathcal{D}(\mathscr{A}_i)$
\end{lemme}
\preuve .
Soit $X\in{\rm Ob}(\mathcal{D}({A}))$. On peut clairement \'ecrire 
$X = \oplus X_i$ où chaque complexe $X_i$ est form\'e d'objets de 
$\mathscr{A}_i$.
Soient maintenant  $X_i$ et $X_j$ des complexes form\'es respectivement 
d'objets de $\mathscr{A}_i$ et $\mathscr{A}_j$.
Soit $I_i$ un complexe form\'e d'injectifs de $\mathscr{A}_i$ 
quasi-isomorphe au complexe $X_i$. 
Par hypoth\`ese, $I_i$ est un complexe form\'e d'injectifs de $\mathscr{A}$, donc 
${\rm{Hom}}_{\mathcal{D}(\mathscr{A})} (X_i,X_j) = 
                          {\rm{Hom}}_{K(\mathscr{A})}(I_i, X_j) = 0$.
On a donc bien une d\'ecomposition de la cat\'egorie d\'eriv\'ee 
$\mathcal{D}({\mathscr{A}})$.
\sq

On applique ce lemme \`a la  cat\'egorie $G{\rm -Qco}(X)$ des $G$-faisceaux 
quasi-coh\'erents sur $X$
, dans le cas où  $G$ agit trivialement sur $X$. 
En effet, 
\begin{proposition}~\label{proposition : D\'ecomposition G-Coh}
Si $G$ agit trivialement sur le sch\'ema $X$, 
on a une d\'ecomposition de  $G{\rm -Qco}(X)$ suivant les repr\'esentations 
irr\'eductibles de $G$, 
dont chaque facteur est isomorphe \`a  ${\rm Qco}(X)$ : Tout $G$-faisceau 
se d\'ecompose en une somme directe 
$\oplus_{\rho \in {\rm Irr}(G)} E_\rho \otimes\rho$, où les $E_\rho$ sont des 
faisceaux quasicoh\'erents sur $X$ (sans action de $G$).
\end{proposition}
\preuve .
Pour commencer, d\'emontrons la d\'ecomposition pour les 
$G$-faisceaux coh\'erents.
Le cas quasi-coh\'erent s'en d\'eduira par limite inductive.
Pour tout $\rho$, repr\'esentation irr\'eductible de $G$, et pour tout  
$G$-faisceau coh\'erent $F$, 
notons 
$$F_\rho = %{\mbox{\PolicePourFophie Hom } }
\mathscr{H}om_{\mathscr{O}_X-G}(\mathscr{O}_X\otimes \rho, F).$$
$F_\rho$ est un $\mathscr{O}_X$-module coh\'erent (sans action de $G$) 
et on a un morphisme naturel
$$\begin{diagram} \oplus_{\rho \in {\rm{Irr}}(G)}(F_\rho\otimes \rho)&\rTo^{\varphi}& F\\
                    f\otimes v&\rTo& f(1\otimes v).
\end{diagram}$$

Montrons que $\varphi$ est un isomorphisme. 

Dans un premier temps, supposons $F$ localement libre. 
Alors le probl\`eme \'etant local, on peut se placer sur un ouvert 
où $F$ est libre et alors pour toute repr\'esentation $\rho$ irr\'eductible de
$G$,  $F_\rho$ est libre.
Sur chaque fibre on a donc un morphisme de 
($G$-$\mathscr{O}_{X,x}$)-module libres 
de type fini qui est un  isomorphisme lorsqu'on passe 
au quotient par l'id\'eal maximal.
En effet, cette d\'ecomposition est connue pour les
repr\'esentation lin\'eaires de  dimensions finies. D'apr\`es le 
lemme de Nakayama, $\varphi $ induit donc un isomorphisme sur chaque fibre. 
Cela implique que $\varphi$ est un isomorphisme.

Maintenant, soit $F$ un $G$-faisceau coh\'erent quelconque et 
soit 
$L_1 \rightarrow L_0\rightarrow F \rightarrow 0$ une pr\'esentation libre 
de $F$.
Si on applique le foncteur 
$\begin{diagram}[size=2em]F &\rTo& \oplus_{\rho \in {\rm{Irr}}(G)}(F_\rho\otimes \rho)\end{diagram}$
\`a cette suite exacte, on obtient une suite exacte car le 
foncteur $\begin{diagram}[size=2em]F&\rTo& F_\rho\end{diagram}$ est exact \`a droite. En effet, 
$ \mathscr{O}_X \otimes \rho$ est un ($\mathscr{O}_X$-$G$)-module projectif 
car il est libre.
Puis, d'apr\`es le lemme des cinq, 
on obtient l'isomorphisme pour  le faisceau $F$, et la proposition est 
d\'emontr\'ee pour les faisceaux coh\'erents..

%Montrons maintenant la d\'ecomposition dans le cas quasi-coh\'erent.
%Cela va être une cons\'equence du fait suivant : $F$ est la limite inductive de 
%ses sous-$G$-faisceaux coh\'erents.
Soit maintenant $F$ un $G$-faisceau quasi-coh\'erent. 
Localement, $F$ est un ($A$-$G$)-module.
Montrons qu'il est limite inductive de ($A$-$G$)-modules 
de types finis.
Comme $G$ est fini, pour tout $x\in F$, $x$ appartient \`a un 
sous-($A$-$G$)-module de type fini, par exemple 
celui engendr\'e par les \'el\'ements de l'orbite de $x$.
On a donc bien : $F = {\rm {limind}}(M)$, pour $M$ 
sous-($A$-$G$)-module de type fini.
On \'ecrit alors $F = {\rm{limind}}(\oplus M_\rho\otimes \rho)$, et comme 
une somme directe est une limite inductive 
sur une cat\'egorie discr\`ete, limite inductive et somme directe commutent, 
de sorte que l'on a 
$F = \oplus_{\rho \in {\rm{Irr}}(G)} {\rm{limind}}(M_\rho)\otimes \rho$ 
qui est bien une d\'ecomposition de la forme attendue.
Le caract\`ere fonctoriel de la construction locale assure que 
les modules ainsi obtenus se recollent pour former des faisceaux de ($\mathscr{O}_X$-$G$)-modules $F_\rho$ tels que 
$$F = \oplus_{\rho \in {\rm{Irr}}(G)} F_\rho\otimes \rho.$$
%\vspace{-0.5 cm}
\sq
\newpage

Les hypoth\`eses du lemme~\ref{dec A->dec D(A)} \'etant v\'erifi\'ees, on a une d\'ecomposition 
de $\mathcal{D}(G{\rm -Qco}(X))$, qui par restriction aux complexes de 
faisceaux \`a homologie coh\'erente donne :

\begin{cor}
Si $G$ agit librement sur le sch\'ema $X$, on a une d\'ecomposition 
$$\mathcal{D}^G(X) = \oplus_{\rho \in {\rm Irr}(G)} \mathcal{D}^\rho(X),$$ 
où pour tout $\rho $, le foncteur 
$\begin{diagram}\mathcal{D}(X)&\rTo^{\otimes \rho}&\mathcal{D}^\rho(X)
\end{diagram}$ est une \'equivalence de cat\'egorie. Ainsi, un objet 
$F\in\mathcal{D}^G(X)$ est image r\'eciproque d'un objet de $\mathcal{D}(X/G)
= \mathcal{D}(X)$
si et seulement si il appartient \`a la composante $\mathcal{D}^{\rho_0}(X)$ 
correspondante \`a la repr\'esentation 
triviale de $G$.
\end{cor}

%%%%%%%%%%%%%%%%%%%%%%%%%%%%%%%%%%%%%%%%%%%%%%%%%%%%%%%%%%%%%%%%%%%%%%%%%%%%%%%

%%%%%%%%%%%%%%%%%%%%%%%%%%%%%%%%%%%%%%%%%%%%%%%%%%%%%%%%%%%%%%%%%%%%%%%%%%%%%%%
\subsection{Description de la cat\'egorie $\mathcal{D}^G(\mathbb{P}^n)$}
% d\'eriv\'ee $G$-\'equivariante de l'espace 
%projectif.}

Soit $\mathbb{P}^n = \mathbb{P}(V) = {\rm Proj}(A)$ où $V$ est un espace 
vectoriel de dimension $(n+1)$
et $A$ est l'alg\`ebre sym\'etrique sur $V$.
On cherche une description de la cat\'egorie d\'eriv\'ee born\'ee $G$-\'equivariante 
$\mathcal{D}(\mathbb{P}^n)$
 des faisceaux coh\'erents sur $\mathbb{P}^n$.

Soit $M_A^G([0,n])$ la cat\'egorie des ($A$-$G$)-modules libres de type fini 
dont les g\'en\'erateurs sont de degr\'es compris entre $0$ et $n$. Soit 
$K_A^G([0,n])$ la cat\'egorie triangul\'ee des complexes
d'objets de $M_A^G([0,n])$ \`a homotopie pr\`es. Comme dans l'article 
\cite{Bej:79} 
de Be{\u\i}linson pour d\'ecrire la cat\'egorie 
d\'eriv\'ee  $\mathcal{D}(\mathbb{P}^n)$  on construit un foncteur 
$\begin{diagram}[size=2em]F : K_A^G([0,n])&\rTo&  \mathcal{D}^G(\mathbb{P}^n)\end{diagram}$.
Dans cette partie, nous adaptons la d\'emarche de Be{\u\i}linson pour d\'emontrer le 
th\'eor\`eme suivant :

\begin{thm}\label{BeilinsonG}
Le foncteur  $\begin{diagram}[size=2em]F : K_A^G([0,n]) &\rTo& \mathcal{D}^G(\mathbb{P}^n)\end{diagram}$ est  
une \'equivalence de cat\'egories.
\end{thm}

%%%%%%%%%%%%%%%%%%%%%%%%%%%%%%%%%%%%%%%%%%%%%%%%%%%%%%%%%%%%%%%%%%%%%%%%%%%%%%%

Pour cela, on va mettre en bijection des objets des deux cat\'egories, 
qui forment des familles g\'en\'eratrices au sens suivant : 
La plus petite sous-cat\'egorie triangul\'ee pleine  de  $K_A^G([0,n])$ 
(respectivement $ \mathcal{D}^G(\mathbb{P}^n)$ ) 
contenant cette classe d'objet est la cat\'egorie enti\`ere.
Rappelons alors le lemme suivant (\cite{Bej:79}):

%%%%%%%%%%%%%%%%%%%%%%%%%%%%%%%%%%%%%%%%%%%%%%%%%%%%%%%%%%%%%%%%%%%%%%%%%%%%%%%

\begin{lemme}\label{eq}
Soient $\mathscr{C}$ et $\mathscr{D}$ deux cat\'egories triangul\'ees, et 
$\begin{diagram}[size=2em]F :\mathscr{C} &\rTo&\mathscr{D }\end{diagram}$ un foncteur triangul\'e.
Si $X_i$ est  une famille g\'en\'eratrice d'objets de $\mathscr{C}$ telle que 
$F(X_i)$ est une famille g\'en\'eratrice de $\mathscr{D}$, si 
pour tout $(i,j)$ et pour tout $l$ : 
 $$\begin{diagram}
F : {\rm Hom}^l(X_i, X_j)&\rTo&{\rm Hom}^l(F(X_i),F(X_j) )\end{diagram}$$ est un 
isomorphisme, alors $F$ est une \'equivalence de  cat\'egorie.
\end{lemme}

%%%%%%%%%%%%%%%%%%%%%%%%%%%%%%%%%%%%%%%%%%%%%%%%%%%%%%%%%%%%%%%%%%%%%%%%%%%%%%%

Pour appliquer ce lemme, il suffit donc de d\'eterminer 
une famille g\'en\'eratrice $\{X_i\}$ de  $K_A^G([0,n])$ et un foncteur 
$\begin{diagram}[size=2em]F : K_A^G([0,n])&\rTo & \mathcal{D}^G(\mathbb{P}^n)\end{diagram} $
v\'erifiant les hypoth\`eses  et tel que 
 $\{F(X_i)\}$  est une famille g\'en\'eratrice de  $\mathcal{D}^G(\mathbb{P}^n)$.
%%%%%%%%%%%%%%%%%%%%%%%%%%%%%%%%%%%%%%%%%%%%%%%%%%%%%%%%%%%%%%%%%%%%%%%%%%%%%%%
%%%%%%%%%%%%%%%%%%%%%%%%%%%%%%%%%%%%%%%%%%%%%%%%%%%%%%%%%%%%%%%%%%%%%%%%%%%%%%%

\subsubsection{Th\'eor\`eme de structure des ($A$-$G$)-modules libres 
de types finis sur une $k$-alg\`ebre gradu\'ee $A$ de type fini}
Dans ce paragraphe,  $A$  est une $k$-alg\`ebre gradu\'ee sur laquelle 
agit un groupe fini $G$. Un ($A$-$G$)-module gradu\'e est un
($A$-$G$)-module gradu\'e en tant que $A$-module et tel que l'action de $G$ 
respecte le degr\'e.

\begin{proposition}\label{llgradu\'e}
Soit $M$ un ($A$-$G$)-module gradu\'e libre de type fini sur $A$
dont les g\'en\'erateurs sont de degr\'es compris entre $0$  et $n$.
Alors il existe des repr\'esentations $V_\alpha\in R(G)$ telles que 
$$M = \oplus _{\alpha =0}^n V_\alpha \otimes A(-\alpha),$$
où $A(-\alpha)$ est l'anneau gradu\'e obtenu en d\'ecalant la graduation de 
$A$ de $\alpha$ : $A(-\alpha)_i  =A_{i-\alpha}$
\end{proposition}

%%%%%%%%%%%%%%%%%%%%%%%%%%%%%%%%%%%%%%%%%%%%%%%%%%%%%%%%%%%%%%%%%%%%%%%%%%%%%%%
\preuve .
Soit $\{e_i, i\in I\}$ une famille  de g\'en\'erateurs de $M$. 
On peut partitionner 
$I = \sqcup_{\alpha\in [O,n]} I_\alpha$ de telle sorte que
pour tout  
$i\in I_\alpha$, le g\'en\'erateur $e_i$ soit de degr\'e $\alpha$.
Pour tout $g$ dans $ G$ et pour tout $\alpha$, si  $i\in I_\alpha$ alors 
$g\cdot e_i$ est de degr\'e $\alpha$, donc : 
$$g\cdot e_i = \sum_{j\in I_\alpha} r_{j,i} e_j + \sum_{j\in I_{\alpha -1}} 
                             r_{j,i} e_j +\ldots +\sum_{j\in I_0} r_{j,i} e_j$$
où $r_{j,i}$ est de degr\'e $\alpha -\deg(e_j)$.
Il s'agit de d\'emontrer que l'on peut choisir la base  $\{e_i\}_{i\in I}$ 
telle que 
pour tout $g$ dans $ G$, pour tout $\alpha$, et pour tout  $i\in I_\alpha$, 
$$g\cdot e_i = \sum_{j\in I_\alpha} r_{j,i} e_j.$$

Soient $\alpha_o = \max\{\alpha \ |\ I_\alpha \neq 0\}$, 
$J = \sqcup_{\alpha <\alpha_o}I_\alpha$ et $N = \oplus_{j\in J}Ae_j$. 
Pour tout $g\in G$,
$g\cdot N\subset N$, de sorte que $N$ est un sous-($A$-$G$)-module de $M$.

On va d\'emontrer qu'il existe un sous-($A$-$G$)-module gradu\'e $N'$ de $M$, 
dont les g\'en\'erateurs sont de degr\'e $\alpha_o$ et  
tel que $M = N\oplus N'$.

Le $A$-module $M$ peut être vu comme $k$-espace vectoriel, car $A$ est une 
$k$-alg\`ebre gradu\'ee ($A_0 = k$).
Comme la caract\'eristique de $k$ ne divise pas l'ordre du groupe  $G$, les 
repr\'esentations de
dimension finie de $G$ sont semi-simples. 
Or  en tant que repr\'esentation de 
$G$, 
$M = \oplus_{n\in \mathbb{N}} M_n$, où $M_n$ est la partie de $M$ de degr\'e 
$n$ : 
$$M_n = \oplus _{\alpha\leq n} \oplus_{i\in I_{\alpha}} A_{n-\alpha}e_i.$$
Comme $M$ est de type fini sur $A$, $I_\alpha$ est fini, et chaque $A_i$ 
\'etant une repr\'esentation de dimension finie de $G$, 
$M_n$ est une repr\'esentation  de dimension finie pour tout $n\in \mathbb{N}$. 
(Remarquons que $M_n$ est stable par $G$ car l'action respecte la graduation).
Ainsi, en tant que somme directe de repr\'esentations semi-simples, $M$ est 
semi-simple.

En particulier, la suite exacte de ($k$-$G$)-modules 
$$\begin{diagram} 0 &\rTo& N&\rTo& M&\rTo& M/N&\rTo& 0\end{diagram}$$ est scind\'ee.
Ainsi, la ($A$-$G$)-base  $\{\bar{e}_i, i\in I_{\alpha_o} \}$ de $M/N$ se 
rel\`eve en $\{e'_i, i\in I_{\alpha_o}\}$
dans $M$, de telle façon que 
$g\cdot e'_i  = \sum_{j\in I_{\alpha_o}} r_{j,i} e'_j$. Le $A$-module
engendr\'e par la famille $\{e'_i, i\in I_{\alpha_o}\}$ est donc un 
sous-($A$-$G$)-module de $M$. 
Cependant, il n'est pas gradu\'e a priori. On a pour tout $i\in I_{\alpha_o}$,
$$e'_i = e_i + \sum_{j\in J} \lambda_{i,j} e_j.$$
Posons  
$$e''_i =[e'_i]_{\alpha_o} = 
e_i + \sum_{\alpha < \alpha_o} \sum_{j\in I_\alpha} 
                              [\lambda_{i,j}]_{\alpha_o -\alpha} e_j.$$
Pour tout $i\in I_{\alpha_o}$, $e''_i$ est un ant\'ec\'edent de $\bar{e}_i\in M/N$ 
de degr\'e $\alpha_o$ et on a :
$$
\begin{array}{rcll}
g\cdot e_j'' &=& g\cdot[ e'_i]_{\alpha _o} \\
     &  =&  [\sum_{j\in I_{\alpha_o}} r_{j,i} e'_j]_{\alpha_o} \\
     & = &  \sum_{j\in I_{\alpha_o}} r_{j,i} [e'_j]_{\alpha_o} &
                              \mbox{(les $r_{j,i}$ sont de degr\'e nul) }\\
     & = &  \sum_{j\in I_{\alpha_o}} r_{j,i} e''_j
\end{array}
$$

Soit $N'$ le $A$-module engendr\'e par la famille 
$\{e''_i, i\in I_{\alpha_o}\}$. $N'$ est un sous-($A$-$G$)-module
de $M$ suppl\'ementaire de $N$. De plus, on a 
$N' = V_{\alpha_o}\otimes A(-\alpha _o)$ où  
$ V_{\alpha_o} =\oplus_{i\in I_{\alpha_o}}k e_i$
est la repr\'esentation de $G$ de dimension finie 
d\'efinie par l'action de $G$ sur la base $\{e''_1, i\in I_{\alpha_o}\}$.  
En effet, $A(-\alpha_o)$ est le module libre de dimension $1$ sur 
$A$ avec un g\'en\'erateur en degr\'e $\alpha_o$.

On peut recommencer le raisonnement en remplaçant
$M$ par $N$. Par r\'ecurrence descendante sur le degr\'e maximal des 
g\'en\'erateurs de $M$, on obtient le r\'esultat :
$$M = \oplus _{\alpha =0}^n V_\alpha \otimes A(-\alpha).$$ \sq

%%%%%%%%%%%%%%%%%%%%%%%%%%%%%%%%%%%%%%%%%%%%%%%%%%%%%%%%%%%%%%%%%%%%%%%%%%%%%%%

%%%%%%%%%%%%%%%%%%%%%%%%%%%%%%%%%%%%%%%%%%%%%%%%%%%%%%%%%%%%%%%%%%%%%%%%%%%%%%%
%%%%%%%%%%%%%%%%%%%%%%%%%%%%%%%%%%%%%%%%%%%%%%%%%%%%%%%%%%%%%%%%%%%%%%%%%%%%%%%

\subsubsection{D\'emonstration du th\'eor\`eme~\ref{BeilinsonG}}

D'apr\`es la proposition~\ref{llgradu\'e}, la famille 
$\{V\otimes A(-i), V\in {\rm Irr}(G), i\in [0,n] \}$ engendre 
la cat\'egorie $M_A^G([0,n])$. Si on la consid\`ere comme une
famille de complexes (où les faisceaux 
sont vus comme des complexes concentr\'es en degr\'e nul), cette famille
engendre la cat\'egorie triangul\'ee $K_A^G([0,n])$.

Consid\'erons le  foncteur additif 
$$\begin{diagram}
F:& M_A^G([0,n]) &\rTo &G{\rm -Coh}(\mathbb{P}^n)\\
&M =\oplus_{\alpha =0}^n V_\alpha \otimes A(-\alpha) & \rTo &
 \oplus_{\alpha =0}^n V_\alpha \otimes \mathscr{O}(-\alpha) \end{diagram}
$$ Il se prolonge \`a la cat\'egorie triangul\'ee des complexes $K_A^G([0,n])$, et 
quitte \`a composer avec le foncteur de localisation, on peut voir $F$ comme \'etant un foncteur de $K_A^G([0,n])$ dans $\mathcal{D}^G(\mathbb{P}^n)$.
D\'emontrons que pour tout triplet $(V, W, l)$, on a 
 
$$
{\rm Hom}^l_{K_A^G([0,n])}(V\otimes A(-i), W\otimes A(-j))\cong 
{\rm Hom}^l_{\mathcal{D}^G(\mathbb{P}^n)}(
V\otimes \mathscr{O}(-i),W\otimes \mathscr{O}(-j)).\quad (1)$$
Soient $(V,W)$ un couple de repr\'esentations irr\'eductibles et 
$(i,j)$ un cou\-ple d'entiers entre $1$ et $n$. Alors d'apr\`es le lemme 
de Schur et  si l'on suppose $k$ alg\'ebriquement clos (de sorte que pour 
toute repr\'esentation irr\'eductible $V$ de $G$, ${\rm{End}}_G(V) = k$)  ,
$${\rm Hom}(V\otimes A(-i), W\otimes A(-j)) = 
\left \{\begin{array}{ll}   0 &\mbox{ si } j > i \\
      0 &\mbox{ si } i = j \mbox{ et } V\neq W\\
      k & \mbox{ si }   i = j \mbox{ et } V = W\\
      k^m &\mbox{ si } j < i \mbox{ et } m(V, A_{i-j}\otimes W) = m,\\
\end{array}\right.$$
où $m(V,E)$ est la multiplicit\'e de la repr\'esentation irr\'eductible $V$ dans la 
repr\'esentation de dimension finie $E$.
Puis, pour $l\neq 0$, comme les morphismes de $K_A([0,n])$ sont des 
morphismes de complexes,
il n'y a pas de morphisme non nul de $V\otimes A(-i)$ vers 
$W\otimes A(-j)$ d\'ecal\'e de $l$ car ils sont concentr\'es en des degr\'es distincts.
$${\rm Hom}^l(V\otimes A(-i), W\otimes A(-j)) = 0.$$
Calculons maintenant le terme de droite dans (1)
Pour $l = 0$, on a 

$$ \begin{array}{rcl}
{\rm Hom}(V\otimes \mathscr{O}(-i), W\otimes \mathscr{O}(-j)) &\cong& 
{\rm Hom}(V\otimes \mathscr{O}_{\mathbb{P}^n}, W\otimes \mathscr{O}(i-j))\\ 
& \cong & {\rm Hom }(V, \Gamma(\mathbb{P}^n, \mathscr{O}(i-j)\otimes W ))\\ 
& \cong &  {\rm Hom }(V, \Gamma(\mathbb{P}^n, \mathscr{O}(i-j) )\otimes W ).\\
\end{array}$$

D'où, d'apr\`es le lemme de Schur,
$${\rm Hom}(V\otimes \mathscr{O}(-i), W\otimes \mathscr{O}(-j)) = 
\left \{\begin{array}{ll}   0 &\mbox{ si } j > i \\
0 &\mbox{ si } i = j \mbox{ et } V\neq W\\
k & \mbox{ si }   i = j \mbox{ et } V = W\\
k^m &\mbox{ si } j < i \mbox{ et } m(V, A_{i-j}\otimes W) = m\\
\end{array}\right.$$

Pour $l\neq 0$, il s'agit de d\'emontrer que 
${\rm Ext}_G^l(V\otimes \mathscr{O}(-i), W\otimes \mathscr{O}(-j))   = 0 $.
Or d'apr\`es \cite{Grothendieck:57}, il s'agit de l'aboutissement de la suite spectrale de terme 
initial
$${\rm E}^{p,q}_2 = 
H^p(G, {\rm Ext}^q(V\otimes \mathscr{O}(-i), W\otimes \mathscr{O}(-j))  ).$$ 
Or comme on est sur un corps $k$ dont la caract\'eristique ne divise pas 
l'ordre de $G$, il n'y a pas de cohomologie des  groupes, donc on a  
$${\rm Ext}_G^l(V\otimes \mathscr{O}(-i), W\otimes \mathscr{O}(-j))   = 
{\rm Ext}^l(V\otimes \mathscr{O}(-i), W\otimes \mathscr{O}(-j))   ^G .$$
Il suffit donc de d\'emontrer que 
${\rm Ext}^l(V\otimes \mathscr{O}(-i), W\otimes \mathscr{O}(-j))   = 0.$
Or il s'agit de  l'aboutissement de la suite spectrale de terme initial 
$${\rm E}^{p,q}_2 = 
H^p(\mathbb{P}^n, %{\mbox{\PolicePourFophie Ext } }
\mathscr{E}xt^q(V\otimes \mathscr{O}(-i), 
                                         W\otimes \mathscr{O}(-j))   ).$$ 
Comme $\mathscr{O}(-i)\otimes V$ est 
localement libre, les Ext locaux sont triviaux pour $q\neq 0$, donc on a 
\\$\begin{array}{rcl}{\rm Ext}^l(V\otimes \mathscr{O}(-i), 
W\otimes \mathscr{O}(-j))  
& = &H^l(\mathbb{P}^n, %{\mbox{\PolicePourFophie Hom } } 
\mathscr{H}om(V\otimes \mathscr{O}(-i), 
W\otimes \mathscr{O}(-j))   )\\ 
& = & H^l (\mathbb{P}^n, V^*\otimes W\otimes \mathscr{O}(i-j)) \\  
& = & H^l (\mathbb{P}^n, \mathscr{O}(i-j))\otimes V^*\otimes W  \\
& = &   0\end{array} $
car $i-j> -n-1$.

Ainsi, pour tout $V,W\in R(G)$ et pour tout $i,j \in [0, n]$, on a 
$${\rm Hom}_{K_A^G([0,n])}(V\otimes A(-i), W\otimes A(-j))\cong 
{\rm Hom}_{\mathcal{D}^G(\mathbb{P}^n)}(F(V\otimes A(-i)),F( W\otimes A(-j))).$$
Pour appliquer le lemme~\ref{eq}, il reste \`a   d\'emontrer 
la proposition suivante :

%%%%%%%%%%%%%%%%%%%%%%%%%%%%%%%%%%%%%%%%%%%%%%%%%%%%%%%%%%%%%%%%%%%%%%%%%%%%%%%
%%%%%%%%%%%%%%%%%%%%%%%%%%%%%%%%%%%%%%%%%%%%%%%%%%%%%%%%%%%%%%%%%%%%%%%%%%%%%%%

\begin{proposition}\label{g\'en\'erationDG(Pn)}
La  famille $ \{V\otimes \mathscr{O}(-i), V\in {\rm Irr}(G), i\in [0,n]  \}$  
engendre la cat\'egorie $\mathcal{D}^G(\mathbb{P}^n)$.
\end{proposition}
Comme dans \cite{Bej:79}, on consid\`ere la r\'esolution de la diagonale par le complexe de 
Koszul :
$$\begin{diagram}
0 &\rTo&p_1^*\Omega^n(n)\otimes p_2^*\mathscr{O}(-n) &\rTo &\ldots\  
 & p_1^*\Omega^1(1)\otimes p_2^*\mathscr{O}(-1)
&\rTo & \mathscr{O}_{\mathbb{P}^n\times \mathbb{P}^n}&\rTo & 
\mathscr{O}_\Delta &\rTo & 0.\end{diagram}$$

La diagonale $\Delta$ est un sous-$G$-sch\'ema de 
$\mathbb{P}^n\times \mathbb{P}^n$, donc $\mathscr{O}_\Delta$ est un 
$G$-faisceau 
et la r\'esolution \'etant compos\'ee de $G$-faisceaux, il s'agit d'une r\'esolution 
dans la cat\'egorie $G{\rm -Coh}(\mathbb{P}^n\times \mathbb{P}^n)$
de $\mathscr{O}_\Delta$ par des $G$-faisceaux  localement libres 
donc projectifs. 
%On montre facilement qu'ils sont \'egalement projectifs 
%dans $G$-$Coh(\mathbb{P}^n\times \mathbb{P}^n)$.
%
%    CI APRES LA PREUVE DE `` G-Faisceau  l.l. => G-faisceau projectif ``: 
%
%Montrons que les $G$-faisceaux localement libres sont acycliques pour le 
%produit tensoriel  $$
%\otimes : G{\rm -Coh}(\mathbb{P}^n\times \mathbb{P}^n)\times 
%G{\rm -Coh}(\mathbb{P}^n\times \mathbb{P}^n)
%\rightarrow G{\rm -Coh}(\mathbb{P}^n\times \mathbb{P}^n).$$
%La d\'emonstration est locale et revient \`a montrer que  tout ($A$-$G$)-modules 
%libre est un ($A$-$G$)-module projectif. 
%Or, il s'agit en particulier d'un $A$-module libre, 
%donc projectif, et tout $A$-module projectif qui est un ($A$-$G$)-module est 
%un ($A$-$G$)-module projectif, car le produit tensoriel de 
%$G{\rm -Coh}(\mathbb{P}^n\times \mathbb{P}^n)$ est une restriction de celui 
%de ${\rm Coh}(\mathbb{P}^n\times \mathbb{P}^n)$.

Consid\'erons les foncteurs 
$$\begin{diagram}
p_1^* \ :\  &G{\rm -Coh}(\mathbb{P}^n) &
\rTo & G{\rm -Coh}(\mathbb{P}^n\times \mathbb{P}^n)
\end{diagram}$$
$$\begin{diagram}
\otimes \mathscr{O}_\Delta\  :\  &G{\rm -Coh}(\mathbb{P}^n\times \mathbb{P}^n)&
\rTo&G{\rm -Coh}(\mathbb{P}^n\times \mathbb{P}^n)
\end{diagram}$$
$$\begin{diagram}
p_{2,*} \ :\  & G{\rm -Coh}(\mathbb{P}^n\times \mathbb{P}^n) &
\rTo & G{\rm -Coh}(\mathbb{P}^n)%\\
\end{diagram}
$$
La composition de ces trois foncteurs est l'identit\'e de 
$G{\rm -Coh}(\mathbb{P}^n)$, de sorte qu'en d\'erivant on obtient :
$$\forall M\in \mathcal{D}^G(\mathbb{P}^n), \  M = 
                  Rp_{2,*}(\mathscr{O}_\Delta\otimes^L p_1^*(M)) .$$
Or, d'apr\`es la r\'esolution de la diagonale par des complexes 
$\otimes$-acycliques, on sait calculer le complexe
$\mathscr{O}_\Delta\otimes^L p_1^*(M)$. En degr\'e $-\alpha$, 
on trouve le faisceau
$\oplus_{d=0}^n p_2^*\mathscr{O}(-d)\otimes p_1^*\Omega^d(d)\otimes 
p_1^*(M)^{d-\alpha}$.
Puis, comme $p_1^*$ est exact et commute au produit tensoriel,   
$  p_1^*\Omega^d(d)\otimes p_1^*(M)^{d-\alpha}= 
p_1^*(\Omega^d(d)\otimes M^{d-\alpha})$. Ainsi, le complexe
$\mathscr{O}_\Delta\otimes^L p_1^*(M)$ appartient \`a la sous-cat\'egorie 
triangul\'ee engendr\'ee par les faisceaux du type :
$p_2^*\mathscr{O}(-d)\otimes p_1^*(N)$, où $N$ est un $G$-faisceau sur 
$\mathbb{P}^n$.

Alors, comme les triangles sont conserv\'es par les  foncteurs d\'eriv\'es, 
on peut obtenir $M$ par une suite de triangles d\'ebutant 
par des objets du type $Rp_{2,*}(p_2^*\mathscr{O}(-d)\otimes p_1^*(N))$, 
où $N\in G{\rm -Coh} (\mathbb{P}^n)$.

Il suffit donc de d\'emontrer que ces objets appartiennent \`a la sous-cat\'egorie 
triangul\'ee engendr\'ee par la famille
$ \{V\otimes \mathscr{O}(-i), V\in {\rm Irr}(G), i\in [0,n]  \}.$
Or, d'apr\`es la formule  de projection, 
$$Rp_{2,*}(p_2^*\mathscr{O}(-d)\otimes p_1^*(N)) = 
\mathscr{O}(-d)\otimes Rp_{2,*}(p_1^*(N)).$$
Par ailleurs, pour tout $i$, $R^ip_{2,*}(p_1^*(N)) = H^i(Rp_{2,*}(p_1^*(N)))$ est le faisceau
 associ\'e au pr\'efaisceau 
$\begin{diagram}[size=2em]V&\rTo &H^i( \mathbb{P}^n \times V, p_1^*(N)_{|\mathbb{P}^n\times V}) = 
H^i(\mathbb{P}^n, P)\end{diagram}$, donc le faisceau
constant $H^i(\mathbb{P}^n, N)$
qui est une repr\'esentation de $G$.

Les faisceaux d'homologie de  
$Rp_{2,*}(p_2^*\mathscr{O}(-d)\otimes p_1^*(N))$ sont de la forme 
$V\otimes \mathscr{O}(-i)$ et appartiennent donc \`a la sous-cat\'egorie 
triangul\'ee voulue.
Le lemme suivant implique qu'alors, le complexe   
$Rp_{2,*}(p_2^*\mathscr{O}(-d)\otimes p_1^*(N))$ appartient aussi 
\`a cette  sous-cat\'egorie.
\sq

\begin{lemme}
Soit $\mathscr{A}$ une cat\'egorie ab\'elienne, $\mathscr{C}$ une classe 
d'objets de $\mathscr{A}$ et $<\mathscr{C}>$ la sous
cat\'egorie triangul\'ee de $\mathcal{D}^b(\mathscr{A})$ engendr\'ee par 
$\mathscr{C}$.
Si $M\in \mathcal{D}^b(\mathscr{A})$ a tous ses faisceaux d'homologie 
dans $<\mathscr{C}>$, alors $M\in <\mathscr{C}>$.
\end{lemme}
\preuve .
On d\'emontre ce lemme par r\'ecurrence sur l'amplitude homologique de $M$.
Si elle est nulle, il n'existe qu'un faisceau  d'homologie non nul 
$H^k(M)$, et $M$ est quasi-isomorphe
au complexe $H^k(M)[-k]$. Puis, comme $<\mathscr{C}>$ est stable par 
triangles, elle l'est en particulier par quasi-isomorphisme, donc 
$M\in <\mathscr{C}>$.

Supposons le lemme d\'emontr\'e pour tout complexe d'amplitude homologique 
inf\'erieure ou \'egale \`a $d$.
Soit $M\in \mathcal{D}^b(\mathscr{A})$ un complexe d'amplitude $d+1$ 
dont les faisceaux d'homologie appartiennent \`a $<\mathscr{C}>$.
Soit $k_o = \min\{k\mbox{ tel que } H^k(M)\neq 0\}$.
Consid\'erons le triangle 
$$\begin{diagram}H^{k_o}(M)&\rTo& M&\rTo& \tau^{>k_o }M\end{diagram}$$
dans $ \mathcal{D}^b(\mathscr{A})$.
Le premier terme appartient \`a  $<\mathscr{C}>$ par hypoth\`ese et le 
troisi\`eme par hypoth\`ese de r\'ecurrence, car $ \tau^{>k_o }M$ est 
d'amplitude homologique $d$ et ses faisceaux d'homologie sont ceux 
de $M$, donc appartiennent \`a $<\mathscr{C}>$ par hypoth\`ese.
Ainsi, comme $<\mathscr{C}>$ est une sous-cat\'egorie triangul\'ee de 
$\mathcal{D}^b(\mathscr{A})$, $M\in <\mathscr{C}>$.
\sq

On peut donc appliquer le lemme~\ref{eq}, et on obtient le 
th\'eor\`eme~\ref{BeilinsonG}
\sq

%%%%%%%%%%%%%%%%%%%%%%%%%%%%%%%%%%%%%%%%%%%%%%%%%%%%%%%%%%%%%%%%%%%%%%%%%%%%%%%

%%%%%%%%%%%%%%%%%%%%%%%%%%%%%%%%%%%%%%%%%%%%%%%%%%%%%%%%%%%%%%%%%%%%%%%%%%%%%%%
%
%        SSS   EEEE   CCC   TTTTTT  I    OOO    N   N
%       S      E     C        T     I   O   O   NN  N
%        SS    EE    C        T     I   O   O   N N N
%          S   E     C        T     I   O   O   N  NN
%       SSSS   EEE    CCC     T     I    OOO    N   N

%%%%%%%%%%%%%%%%%%%%%%%%%%%%%%%%%%%%%%%%%%%%%%%%%%%%%%%%%%%%%%%%%%%%%%%%%%%%%%%

\section{Crit\`ere de descente d'un objet de la cat\'egorie d\'eriv\'ee $G$-\'equivariante }

Dans cette deuxi\`eme partie, nous \'evaluons la distance qui s\'epare les cat\'egories
$\mathcal{D}^G(X)$ et $\mathcal{D}(X/G)$. Autrement dit, nous d\'eterminons 
l'obstruction pour un objet de $\mathcal{D}^G(X)$ \`a \^etre l'image par  $\pi^*$ 
d'un objet de $\mathcal{D}(X/G)$

On a vu dans le th\'eor\`eme~\ref{equivalence_cat_der_sous_cat_loclibre} 
que $\mathcal{D}^G(X)\simeq \mathscr{M}$, mais $\mathscr{M}$ 
n'est pas la cat\'egorie d\'eriv\'ee born\'ee de la cat\'egorie $\mathscr{L}$
des $G$-faisceaux coh\'erents localement libres sur $X$.

Dans un premier paragraphe, on traite malgr\'e cela le probl\`eme 
sur $\mathscr{L}$. 

On se ram\`enera alors \`a ce cas  
en d\'emontrant qu'un objet de $\mathscr{M}$ est image par $\pi^*$ d'un
objet de $\mathcal{D}(X/G)$ si et seulement si chacun des faisceaux de $\mathscr{L}$ 
qui le compose est image par$\pi^*$ d'un faisceau de
${\rm Coh}(X/G)$.
Cela r\'esout le probl\`eme de la descente d'un objet de $\mathcal{D}^G(X)$.

\subsection{Descente d'un $G$-faisceau localement libre.}

On commence par raisonner dans le cas local. 
La proposition suivante traite le cas de la descente d'un $G$-faisceau
localement libre au voisinage d'un point  ferm\'e  de $X$ (sch\'ema projectif
lisse sur un corps $k$ alg\'ebriquement clos de caract\'eristique nulle).
\begin{proposition}\label{R-G_module_libre}
Soit $R$ une $k$-alg\`ebre locale d'id\'eal maximal $\mathfrak{m}$ et de corps r\'esiduel $k$.
Alors la cat\'egorie $\mathscr{L}$ des ($R$-$G$)-modules libres est 
\'equivalente \`a la  cat\'egorie $\mathscr{R}$ des repr\'esentations de $G$.
Les foncteurs d'\'equivalence sont :
$$\begin{diagram}
\mathscr{L}&\rTo^{\otimes k} &\mathscr{R} \end{diagram}{\rm{\ et \ }}
\begin{diagram}
\mathscr{R} &\rTo^{\otimes R} &\mathscr{L}.\end{diagram}$$
\end{proposition}

\preuve.
On a clairement 
$$\forall V\in\mathscr{R}, \ (V\otimes_k R)\otimes_R k = V\otimes_k R/\mathfrak{m}\cong V.$$
R\'eciproquement, il s'agit de d\'emontrer que 
$$\forall M\in\mathscr{L},\   M/\mathfrak{m}M \otimes_k R\cong M {\mbox{ dans }}\mathscr{L}.$$

Soit donc $M$ un ($R$-$G$)-module libre. Comme $M$ est 
la limite inductive de ses sous-($R$-$G$)-modules de types  finis, on peut 
supposer $M$ de type fini.
En effet, si $M = {\rm limind}(N)$, alors 
$$(M\otimes k) \otimes R =  ({\rm limind }(N )\otimes k)\otimes R = 
{\rm{limind}}(N\otimes k)\otimes R = {\rm limind}((N\otimes k)\otimes R)$$
car les limites inductives commutent  aux foncteurs adjoints \`a gauches, 
donc en particulier aux produits tensoriels.

Soit $V= M\otimes k = M/\mathfrak{m}M$. On veut d\'emontrer que $M\cong V\otimes R$
en tant que ($R$-$G$)-module.
$V$ \'etant une repr\'esentation de type fini de $G$, on peut en choisir une base 
$\overline{e_1},\ldots, \overline{e_n}$, où les $e_i$ sont dans $M$ et 
$\overline{e_i}$ est la classe modulo $\mathfrak{m}$. Soit $W$ le 
sous-($k$-$G$)-module de $M$ engendr\'e par les $e_i$. C'est le sous-espace vectoriel
engendr\'e par les $e_i$ et leurs transform\'es par $G$. 
Il est  de dimension finie car $G$ est fini. 
Soit $V_1$ le quotient de $W$ par un suppl\'ementaire de $V$ en tant que $k$-espace vectoriel. 
La surjection  $\begin{diagram}[size=2em]W &\rTo& V_1\cong V\end{diagram}$,  
est un morphisme  de ($k$-$G$)-modules par construction. 
Soit $V'$ son noyau. On a  $$\begin{diagram}0 &\rTo & V'&\rTo & W &\rTo &V &\rTo &0\end{diagram}.$$
Comme toute repr\'esentation de  
dimension finie, $W$ est  semi-simple. 
Il existe donc un ($k$-$G$)-module 
$V_2\cong V$ telle que $W = V_2\oplus V'$.
Notons le encore $V $.
On a alors un morphisme naturel de ($R$-$G$)-modules
$$\begin{diagram}V\otimes_k R&\rTo^{\varphi} &M\end{diagram}.$$ 
 
D'apr\`es le lemme de Nakayama, il est surjectif, et comme il 
s'agit d'une surjection entre $R$-modules libres de même rang, c'est un isomorphisme.
\sq

\begin{cor}\label{isol\'e}
Soit $E$ un $G$-faisceau localement libre sur $X$. 
Il existe un faisceau  $F\in{\rm{Qco}} (X/G)$ tel que  $E=\pi^*(F)$ 
si et seulement si en chacun des points fixes ferm\'es de $X$, 
la repr\'esentation $i^*_x(E) = E_x/\mathfrak{m}_x E_x$ du stabilisateur 
de $x$ est  triviale.
\end{cor}

\preuve. Le morphisme canonique $\begin{diagram}[size=2em]\pi^*\pi_*^G(E)&\rTo & E\end{diagram}$ est un 
isomorphisme si et seulement si
c'est un isomorphisme sur chacune des fibres
%Il suffit en fait de tester les fibres au dessus des points ferm\'es par 
%densit\'e .
au dessus des points ferm\'es.

En tout point où $G$ agit librement, il s'agit bien d'un isomorphisme 
d'apr\`es le  paragraphe~\ref{actionlibre}.

Soit  $x$ est un point fixe de $X$ et $H$ son stabilisateur.
D'apr\`es la proposition~\ref{R-G_module_libre}, si $x$ est un point ferm\'e on a  
$$E_x \cong \mathscr{O}_{X,x}\otimes V$$ où $V = E_x/\mathfrak{m}_xE_x.$
Par ailleurs,  
$\pi_*^G(E) = E_x^H$ est un $\mathscr{O}_{X/G,\pi(x)}$-module libre
donc d'apr\`es  la proposition~\ref{R-G_module_libre}, 
$$E_x^H \cong  \mathscr{O}_{X/G, \pi(x)}\otimes E_x^H/\mathfrak{m}_{\pi(x)} E_x^H 
        \cong  \mathscr{O}_{X/G, \pi(x)}\otimes V^H.$$
Si $\mathscr{O}_{X,x}$ est complet, on a $ \mathscr{O}_{X, \pi(x)} = \mathscr{O}_{X,x}^H$,
donc on obtient $E_x^H =\mathscr{O}_{X,x}^H\otimes V^H$, d'où
$$\begin{array}{rcl}
(\pi^*\pi_*^G(E))_x  &\cong &
( \mathscr{O}_{X,x}^H \otimes V^H) \otimes_{\mathscr{O}_{X,x}^H}\mathscr{O}_{X,x}\\
                     &\cong&  V^H \otimes \mathscr{O}_{X,x} 
\end{array}$$
Ainsi, le morphisme induit un isomorphisme sur la fibre en $x$ 
si et seulement si $V^H = V$, autrement dit, 
$i^*_x(E)$ est la repr\'esentation triviale de $H$. 

Sinon, on se ram\`ene \`a ce cas l\`a car $\widehat{\mathscr{O}_ {X,x}}$ est fid\`element plat sur
$\mathscr{O}_ {X,x}$.
\sq

%%%%%%%%%%%%%%%%%%%%%%%%%%%%%%%%%%%%%%%%%%%%%%%%%%%%%%%%%%%%%%%%%%%%%%%%%%%%%%%

%%%%%%%%%%%%%%%%%%%%%%%%%%%%%%%%%%%%%%%%%%%%%%%%%%%%%%%%%%%%%%%%%%%%%%%%%%%%%%
Dans le cas g\'en\'eral, les points fixes ne sont pas isol\'es.
Une condition du type du corollaire~\ref{isol\'e} devient donc 
insuffisante, car trop contraignante : il s'agit de 
v\'erifier une condition sur chaque sous-sch\'ema ferm\'e d\'efini par un 
point fixe de $X$.

On construit une stratification de 
$Fix(X) = \{ x\in X \ |\ G_x \neq 0\}$ de la façon suivante.
Pour tout sous-groupe $H$ de $G$, notons $\Delta_H$ l'ensemble des 
points de stabilisateur $H$ et $X^H$ le sch\'ema des points 
fixes par $H$ :  
$$\Delta_H = \{ x\in X \ |\   G_x = H\}, \ X^H = \{x\in X\ | \ Hx = x \} = 
\bigsqcup_{H\subset K\subset G}\Delta_K.$$
$X^H$ est lisse car il s'agit du sous-sch\'ema des points fixes sous 
l'action du sous-groupe $H$ de $G$ r\'eductif.
%Ainsi, $\Delta_H$ qui est un ouvert de $X^H$ 
%(en effet $\Delta_H = X^H\setminus \bigsqcup_{K\supseteq H\\K\neq H}X^K$)
%
%
%  DDD  EEEE  M   M  OOO   ??
%  D  D E     MM MM O   O    ?
%  D  D EE    M M M O   O  ??
%  D  D E     M   M O   O   
%  DDD  EEEE  M   M  OOO   ?
%

On consid\`ere alors le sous-sch\'ema  ferm\'e $\overline{\Delta}_H$ de $X$ 
muni de sa structure r\'eduite, et on note 
$\begin{diagram}[size=2em]i_H :\overline{\Delta}_H &\rInto &X \end{diagram}$ 
l'inclusion ferm\'ee.
\begin{lemme}
$\overline{\Delta}_H$ est lisse.
\end{lemme}
\preuve . Notons $X^H = \sqcup_{\alpha\in I}X_\alpha^H$ les 
composantes irr\'eductibles de $X^H$. Alors % ici IRREDUCTUBLE = CONNEXE CAR X^H LISSE 
$\Delta_H = \sqcup_{\alpha\in I}\Delta_H \cap X_\alpha^H$. 
Puis, si $J\subset I$ est le sous-ensemble 
(non vide si $ \Delta_H $ est non vide ) 
form\'e des indices $\alpha$ tels que $\Delta_H \cap X_\alpha^H$ est non vide, 
alors $\Delta_H = \sqcup_{\alpha\in J}\Delta_H \cap X_\alpha^H$
 et $\overline{\Delta}_H = \sqcup_{\alpha\in J}X_\alpha^H$. 
En effet, l'adh\'erence d'un sous-ensemble non vide 
d'une composante irr\'eductible $X_\alpha^H$ est $X_\alpha^H$ entier.
Ainsi, $\overline{\Delta}_H$ est la r\'eunion de certaines composantes % POURQUOI PAS DIRE IRR\'EDUCTIBLES ->matthieu
connexes de $X^H$ qui est lisse, donc est lisse.
\sq

%Le groupe $N_G(H)$ agit sur $\overline{\Delta}_H$  de telle façon que  
%$N_G(H)/H$ agit librement sur $\overline{\Delta}_H$, et
 $H$ agit trivialement sur $\overline{\Delta}_H\subset X_H$.
Ainsi, d'apr\`es le paragraphe~\ref{decomposition} on a une d\'ecomposition : 
$$H{\rm -Qco}(\overline{\Delta}_H) = 
\oplus_{\rho\in R(H)} {\rm Qco}^{\rho}(\overline{\Delta}_H)$$

\begin{proposition}\label{faisceaull}
Soit $E$ un $G$-faisceau localement libre sur $X$. 
Pour tout $H$ sous-groupe de $G$,
$i_H^*(E) = E_H$ est un  $H$-faisceau  
localement libre sur $\overline{\Delta}_H$.
Alors $E$ est l'image r\'eciproque d'un faisceau localement libre sur $X/G$ 
si et seulement si pour tout sous-groupe $H$ de $G$,
$i_H^*(E)$ appartient \`a  la  composante correspondant \`a la repr\'esentation 
triviale de $H$ dans la d\'ecomposition de $H {\rm -Qco}(\overline{\Delta}_H)$.
\end{proposition}
\preuve . Soit $H$ un sous-groupe de $G$.
Si $\Delta_H$ est vide, on conviendra que $i_H^*(E)$ est 
trivial dans le sens ci-dessus. Sinon, soit $x\in \Delta_H$.
En particulier, $x$ est un point fixe.
D'apr\`es le corollaire~\ref{isol\'e}, $E$ provient du bas si et seulement 
si pour tout $x\in Fix(X)$, $i_x(E)$ est une repr\'esentation
triviale du stabilisateur de $x$. 
Consid\'erons le morphisme  
$\begin{diagram}[size=2em] i_x^H :
{\rm Spec}(k(x)) &\rTo &\overline{\Delta}_H\end{diagram}$ d'image $x$. 
Alors $i_x = i_H\circ i^H_x$, de sorte que si 
$i_H^*(E)$   appartient \`a  la  composante correspondant \`a la 
repr\'esentation triviale de $H$ dans la 
d\'ecomposition de $H {\rm -Qco}(\overline{\Delta}_H)$, alors
$i^*_x(E) = i_x^{H*}\circ i_H^*(E)$ est une repr\'esentation triviale de $H$.
D'où la condition suffisante.

R\'eciproquement, supposons que $E = \pi^*(F)$ où $F$ est un 
faisceau localement libre sur $X/G$.
$E$ est un $G$-faisceau trivial sur $X$ en ce sens que l'action 
de $G$ sur  $E$ est diagonale : 
Pour tout ouvert   $U\subset X$, le $(G$-$\mathscr{O}_X(U))$-module 
$E(U)$ v\'erifie : 
$g.(ae) = (g.a) e$ pour tout $(g,a,e)\in G\times \mathscr{O}_X (U)\times E(U)$.
Ainsi,  comme l'action de $H$ est triviale sur $\overline{\Delta}_H$, 
elle l'est \'egalement sur $i_H^*(E)$, de sorte que 
$i_H^*(E)$
appartient \`a  la  composante correspondant \`a la repr\'esentation 
triviale de $H$ dans la d\'ecomposition de $H {\rm -Qco}(\overline{\Delta}_H)$.
\sq
 
% car il s'agit de la tensorisation du  faisceau sans action par la 
% repr\'esentation triviale
%%%%%%%%%%%%%%%%%%%%%%%%%%%%%%%%%%%%%%%%%%%%%%%%%%%%%%%%%%%%%%%%%%%%%%%%%%%%%%%

% R\'EDACTION ET ARGUMENT    A REVOIR 

%  DDD  EEEE  M   M  OOO   ??
%  D  D E     MM MM O   O    ?
%  D  D EE    M M M O   O  ??
%  D  D E     M   M O   O   
%  DDD  EEEE  M   M  OOO   ?

%%%%%%%%%%%%%%%%%%%%%%%%%%%%%%%%%%%%%%%%%%%%%%%%%%%%%%%%%%%%%%%%%%%%%%%%%%%%%%%
%%%%%%%%%%%%%%%%%%%%%%%%%%%%%%%%%%%%%%%%%%%%%%%%%%%%%%%%%%%%%%%%%%%%%%%%%%%%%%%

\subsection{R\'eduction du probl\`eme}\label{red2}
On va maintenant d\'emontrer que pour qu'un complexe de faisceaux 
coh\'erents localement libres provienne du bas, il 
suffit que chacun des faisceaux qui le compose soit image r\'eciproque 
d'un faisceau de ${\rm Coh}(X/G)$.

Commençons par d\'emontrer le lemme suivant :

%%%%%%%%%%%%%%%%%%%%%%%%%%%%%%%%%%%%%%%%%%%%%%%%%%%%%%%%%%%%%%%%%%%%%%%%%%%%%%%

\begin{lemme}\label{lemme}

Soit $A$ un anneau local avec action de $G$. Soient $V$ et $W$ deux 
repr\'esentations triviales de dimensions finies du groupe $G$ et 
$\begin{diagram}[size=2em]
V\otimes A &\rTo^\phi& W\otimes A
\end{diagram}$ 
un morphisme de $A$-modules commutant \`a l'action de $G$. 
Alors $\phi$ provient d'un morphisme 
$\begin{diagram}[size=2em] V\otimes A^G &\rTo& W\otimes A^G \end{diagram}$. 
\end{lemme}

\preuve.
Soit $(e_1, \ldots, e_N)$ une base de $V$ et $(f_1, \ldots, f_M)$ une 
base de $W$.
On a   : \\$\phi(\sum_{i=1}^N  a_i e_i) = 
\sum_{i=1}^N a_i\phi(e_i) = \sum_{i=1}^N a_i\sum_{i=1}^M \phi_{j,i} f_j
$
où l'application $A$-lin\'e\-aire $\phi$ est repr\'esent\'ee par la 
matrice $((\phi_{j,i}))\in \mathscr{M}_{M,N} (A)$.

Alors le fait que $\phi $ commute \`a l'action de $G$ s'\'ecrit de la 
mani\`ere suivante 
(compte tenu du fait que $G$ agit trivialement sur $V$ et $W$, et 
donc que les bases sont invariantes)
$$
\begin{array}{lrcll}
\forall g\in G, \quad &
g\cdot\phi(\sum_{i=1}^N  a_i e_i)  &= &
\sum_{i=1}^N g\cdot a_i\sum_{i=1}^M g\cdot\phi_{j,i} f_j& =  \\
&
\phi(g\cdot \sum_{i=1}^N  a_i e_i) & =&
\sum_{i=1}^N g\cdot a_i\sum_{i=1}^M \phi_{j,i} f_j.
\end{array}$$                  
De sorte que $\phi$ commute \`a l'action de $G$ revient \`a dire que 
$\phi$ est $G$-invariant, en ce sens qu'il provient d'un morphisme 
$\begin{diagram}[size=2em] V\otimes A^G &\rTo& W\otimes A^G \end{diagram}$.
\sq

%%%%%%%%%%%%%%%%%%%%%%%%%%%%%%%%%%%%%%%%%%%%%%%%%%%%%%%%%%%%%%%%%%%%%%%%%%%%%%%

\begin{proposition}\label{morphisme}
Soit 
$\begin{diagram}[size=2em] L_1 = \pi^*(E_1)&\rTo^{\phi}& L_2 = 
\pi^*(E_2)\end{diagram}$ un $G$-morphisme entre $G$-faisceaux localement 
libres de types finis sur $X$ 
images r\'eciproques de faisceaux sur $X/G$. Alors $\phi$  provient d'un 
morphisme unique $\begin{diagram}[size=2em] E_1&\rTo& E_2\end{diagram}$ dans 
${\rm{Coh}}(X/G)$.
\end{proposition}
\preuve.
L'identit\'e $\pi_*^G\circ \pi^* = {\rm Id}_{{\rm{Coh}}(X/G)}$
entraine que s'il existe, ce morphisme est n\'ecessairement 
$\begin{diagram}[size=2em]\pi_*^G(L_1) =
E_1 & \rTo^{ \pi_*^G(\phi)}&  \pi_*^G(L_2) = E_2\end{diagram}$. 
Par aileurs, l'adjonction de  $\pi^*$ et $\pi_*^G$ fournit un morphisme 
canonique  
$\begin{diagram}[size=2em]\pi^* \circ \pi_*^G& \rTo^{\eta} &
{\rm Id}_{G-{\rm Coh }(X)}\end{diagram}$ d'où le diagrame :

$$
\begin{diagram}
L_1 &\rTo^{\pi^* (\pi_*^G(\phi))}& L_2\\
\dTo^{\eta (L_1)} & &\dTo_{\eta(L_2)}\\
L_1 & \rTo^{\phi}& L_2
\end{diagram}
$$
Il s'agit donc de d\'emontrer que  $\eta (L_i) = {\rm Id}_{L_i}$, 
pour $i = 1,2$ et que $\phi = \pi^* (\pi_*^G(\phi))$.

On peut raisonner localement, et donc se placer sur un ouvert affine 
${\rm Spec}(A)$, où $A$ est local.
Comme $L_1$ et $L_2$ sont localement libre, la 
proposition~\ref{R-G_module_libre} 
s'applique et donc $L_i = V_i\otimes A$, où $V_i$ est 
une repr\'esentation de $G$. Alors l'hypoth\`ese que les $L_i$ proviennent 
du bas se traduit par le fait que 
ces repr\'esentations sont triviales.
Par ailleurs, $\phi$ \'etant un morphisme de $G$-faisceaux, il commute \`a 
l'action de $G$ et on peut appliquer le lemme~\ref{lemme}.
Ainsi, $\phi$ est $G$ invariante.

Or dans le cas affine, 
$\begin{diagram}[size=2em]{\rm Spec}(A)&\rTo^{\pi} & {\rm Spec}(B)\end{diagram}$, 
où $B = A^G$, l'op\'eration
$\pi_*^G$ revient \`a consid\'erer un $A$-module comme $B$-module et 
prendre sa partie $G$ invariante, de sorte que
le complexe 
$\begin{diagram}[size=2em]  V_1\otimes A&\rTo^{\phi}& V_2\otimes A\end{diagram}$ 
est transform\'e en 
$\begin{diagram}[size=2em]  V_1\otimes A^G&\rTo^{\phi^G = \phi}&
 V_2\otimes A^G\end{diagram}$, puis $\pi^*$ revient \`a tensoriser par 
$A$ au dessus de $A^G$, 
de telle sorte qu'on retrouve le complexe  
$\begin{diagram}[size=2em]  V_1\otimes A&\rTo^{\phi}& V_2\otimes A\end{diagram}$.

On a donc $\pi^*\pi_*^G\phi = \phi$. Autrement dit, $\phi$ provient 
du morphisme $\begin{diagram}[size=2em]\pi_*^G\phi : E_1 &\rTo& E_2\end{diagram}$.
\sq

\begin{cor}
Soit $L$  un objet de $\mathcal{D}^G(X)$ form\'e de $G$-faisceaux 
localement libres.
Alors  $L$ est image r\'eciproque d'un objet de $\mathcal{D}(X/G)$ 
si et seulement si pour tout $i$, $L_i$ est
image r\'eciproque d'un faisceau de $\rm{Coh}(X/G)$.
\end{cor}
\preuve.
La partie directe est \'evidente. R\'eciproquement, si  $L$ est un 
complexe du type : 
$$\begin{diagram} \pi^*(E_1)&\rTo^{\phi_1}&\pi^*(E_2) &\ldots &
\rTo^{\phi_{r-1}} &\pi^*(E_r) &\rTo^{\phi_r}\end{diagram}$$
où chaque $\phi_i$ est un morphisme de $G$-faisceaux, alors en 
appliquant la proposition~\ref{morphisme} on obtient que chaque 
morphisme  provient du bas, 
et donc le complexe provient du bas.
\sq

%%%%%%%%%%%%%%%%%%%%%%%%%%%%%%%%%%%%%%%%%%%%%%%%%%%%%%%%%%%%%%%%%%%%%%%%%%%%%%%
%%%%%%%%%%%%%%%%%%%%%%%%%%%%%%%%%%%%%%%%%%%%%%%%%%%%%%%%%%%%%%%%%%%%%%%%%%%%%%%
%%%%%%%%%%%%%%%%%%%%%%%%%%%%%%%%%%%%%%%%%%%%%%%%%%%%%%%%%%%%%%%%%%%%%%%%%%%%%%%

\subsection{Crit\`ere  g\'en\'eral}
D'apr\`es le paragraphe~\ref{decomposition}, 
pour tout sous-groupe $H$ de $G$ et tout 
sch\'ema $\Delta_H$ associ\'e, on a  une d\'ecomposition de la cat\'egorie 
d\'eriv\'ee $\mathcal{D}^{H}(\overline{\Delta}_H)$ suivant les repr\'esentations
irr\'eductibles de $H$.

\begin{thm}\label{thm_critere_general}
Un objet $E\in\mathcal{D}^G(X)$ est l'image par $\pi^*$  d'un objet de
$\mathcal{D}(X/G)$
si pour tout sous-groupe $H$ de $G$, 
$Li_H^*(E)\in \mathcal{D}^{H}(\overline{\Delta}_H)$ appartient \`a la 
composante correspondant  \`a la repr\'esentation 
triviale dans la d\'ecomposition de la cat\'egorie d\'eriv\'ee. 
Autrement dit,  la fl\`eche 
$$\begin{diagram}\mathcal{D}^G(X)& \rTo^{\oplus_H Li_H^*}  &
\oplus_{H}\mathcal{D}^{H}(\overline{\Delta}_H)\end{diagram}$$ 
est \`a but dans $\oplus_{H}\mathcal{D}^{\rho_0}(\overline{\Delta}_H)$.
\end{thm}
\preuve.
D'apr\`es le  paragraphe~\ref{red1} on peut toujours 
supposer que $E$ est form\'e de $G$-faisceaux coh\'erents localement libres.
Puis, d'apr\`es le paragraphe~\ref{red2} un complexe de faisceaux 
localement libre se descend si et seulement si chacun des faisceaux 
se descend. 
Alors en appliquant la proposition~\ref{faisceaull}, on obtient que   
$E$ provient du bas si et seulement si pour tout sous-groupe $H$ de $G$ et 
pour tout $k\in \mathbb{Z}$, 
$i_H^*(E_k)$ est un $H$-faisceau coh\'erent sur $\overline{\Delta}_H$ 
appartenant \`a la composante correspondant \`a la repr\'esentation triviale dans la
d\'ecomposition de $H{\rm -Coh}(\overline{\Delta}_H)$.
Or les faisceaux localement libres sont $i_H^*$-acycliques, 
%de  sorte que le complexe $i_H^*(E_k), k\in\mathbb{Z}$ obtenu en calculant 
%terme \`a terme est le complexe $Li^*_H(E)$.
donc cette condition s'\'ecrit encore (d'apr\`es le paragraphe~\ref{decomposition})
  $Li^*_H(E) \in \mathcal{D}^{\rho_0}(\overline{\Delta}_H)$.

\sq

%%%%%%%%%%%%%%%%%%%%%%%%%%%%%%%%%%%%%%%%%%%%%%%%%%%%%%%%%%%%%%%%%%%%%%%%%%%%%%%
%%%%%%%%%%%%%%%%%%%%%%%%%%%%%%%%%%%%%%%%%%%%%%%%%%%%%%%%%%%%%%%%%%%%%%%%%%%%%%%

\vspace{1.5 cm}
\noindent
{\bf{Sophie T\'erouanne}}\\
Universit\'e Grenoble 1\\
Institut Fourier   \\           
BP 74               \\       
38402 Saint Martin d'H\`eres Cedex  \\
France\\
Sophie.Terouanne@ujf-grenoble.fr  

\end{document}